\documentclass[11pt]{article}
\usepackage{amsfonts}
\usepackage{latexsym,amsmath}
\usepackage{amssymb,array}
\makeatletter \oddsidemargin 0in \evensidemargin 0in \textwidth
16cm \RequirePackage[dvips]{graphicx} \textheight 20cm
\newtheorem{t1}{Theorem}[section]
\newtheorem{p1}{Proposition}[section]

\newtheorem{c1}{Corollary}[section]
\newtheorem{d1}{Definition}[section]
\newtheorem{r1}{Remark}[section]

\newtheorem{ex}{Example}[section]

\begin{document}
\title{{\bf On Weighted Measure of Inaccuracy for Doubly Truncated Random Variables}}
\author{{\bf Chanchal Kundu}\footnote{Department of Mathematics, Rajiv Gandhi Institute of Petroleum Technology, Rae Bareli 229 316, U.P., India. E-mail: ckundu@rgipt.ac.in; chanchal$_{-}$kundu@yahoo.com.}}
\date{Revised version to appear in {\it Communications in Statistics$-$Theory \& Methods},\\$\copyright$ by The Taylor $\&$ Francis Group.\\ Submitted: June, 2014}
\maketitle

\begin{abstract}
Recently, authors have studied weighted version of Kerridge
inaccuracy measure for truncated distributions. In the present
communication we introduce the notion of weighted interval
inaccuracy measure for two-sided truncated random variables. In
reliability theory and survival analysis, this measure may help to
study the various characteristics of a system/component when it
fails between two time points. Various aspects of weighted
interval inaccuracy measure have been discussed and some
characterization results have been provided. This new measure is a
generalization of recent dynamic weighted inaccuracy measure.
\end{abstract}
\section{Introduction}
The idea of information theoretic entropy was introduced by
Shannon (1948) and Weiner (1949). Shannon was the one who formally
introduced entropy, known as {\it Shannon's entropy} or {\it
Shannon's information measure}, into information theory, and
characterized the properties of information sources and of
communication channels to analyze the outputs of these sources.\\
\hspace*{.2in} Let us consider an absolutely continuous
nonnegative random variable $X$ with probability density function
$f$, distribution function $F$ and survival function $\overline
F\equiv 1-F$. Then the Shannon's information measure or the
differential entropy of $X$ is given by
\begin{eqnarray}\label{eq1}H_X=-\int_0^\infty
f(x)\ln f(x)dx,\end{eqnarray} which measures the expected
uncertainty contained in $f(\cdot)$ about the predictability of an
outcome of $X$.\\
\hspace*{.2in} Since the pioneering contributions by Shannon and
Weiner, numerous efforts have been made to enrich and extend the
underlying information theory. One important development in this
direction is inaccuracy measure due to Kerridge (1961) which can
be thought of as a generalization of Shannon's entropy. It has
been extensively used as a useful tool for measurement of error in
experimental results. Suppose that an experimenter states the probabilities
of the various possible outcomes of an experiment. His statement can lack precision
in two ways: he may not have enough information and so his statement is vague,
or some of the information he has may be incorrect. All statistical inference related
 problems are concerned with making statements which may be inaccurate
in either or both of these ways. Kerridge (1961) proposed the {\it inaccuracy measure} that can
take accounts for these two types of errors. Suppose that the
experimenter asserts that the probability of the $i^{th}$
eventuality is $q_i$ when the true probability is $p_i$. Then the
inaccuracy of the observer can be measured by
$$I(P,Q)=-\sum_{i=1}^n p_i\ln q_i,$$ where
$P=(p_1,p_2,\ldots,p_n)$ and $Q=(q_1,q_2,\ldots,q_n)$ are two
discrete probability distributions such that $p_i\geq0,~q_i\geq0$
and $\sum_{i=1}^np_i=1=\sum_{i=1}^nq_i$.\\
\hspace*{.2in} Nath (1968) extended Kerridge's inaccuracy measure
to the case of continuous situation and discussed some properties.
If $F(x)$ is the actual distribution corresponding to the
observations and $G(x)$ is the distribution assigned by
the experimenter and $f,~g$ are the corresponding density
functions, then the inaccuracy measure is defined as
\begin{eqnarray}\label{eq2}H_{X,Y}=-\int_0^\infty
f(x)\ln g(x)dx.\end{eqnarray} It has applications in statistical
inference and coding theory. When $g(x)=f(x)$, then
(\ref{eq2}) becomes (\ref{eq1}), the Shannon's entropy. The
definition of inaccuracy measure was also extended to truncated
situation, see, Nair and Gupta (2007), Taneja et al. (2009) and Kumar et al.
(2011) for further details.\\
\hspace*{.2in} It is well-known that Shannon entropy is a shift independent measure.
However, in certain applied contexts, such as reliability or mathematical neurobiology,
it is desirable to deal with shift-dependent information measures. Indeed, knowing that
a device fails to operate, or a neuron to release spikes in a given time-interval, yields
relevantly different information from the case when such an event occurs in a different
equally wide interval. In some cases we are thus led to resort to a shift-dependent information
measure that, for instance, assigns different measures to such distributions. Also, there exist
many fields dealing with random experiment whose elementary events are characterized both by
their objective probabilities and by some qualitative (objective or subjective) weights attached
to elementary events and which may or may not be dependent on the objective probabilities.\\
\hspace*{.2in} In analogy with Belis and Guia\c{s}u (1968), Di Crescenzo and Longobardi (2006) considered
the notion of weighted entropy
\begin{eqnarray}\label{eq4}H_{X}^w=-\int_0^\infty
x f(x)\ln f(x)dx.\end{eqnarray}
As pointed out by Belis and Guia\c{s}u (1968) that the occurrence of an event removes a double uncertainty: the quantitative one, related to the probability with which it occurs, and the qualitative one, related to its utility for the attainment of the goal or to its significance with respect to a given qualitative characteristic. The factor $x$, in the integral on the right-hand-side of (\ref{eq4}), may be viewed as a weight linearly emphasizing the occurrence of the event $\{X=x\}$. This yields a length biased shift-dependent information measure assigning greater importance to larger values of $X$. The use of weighted entropy (\ref{eq4}) is also motivated by the need, arising in various communication and transmission problems, of expressing the usefulness of events by means of an information measure.\\
\hspace*{.2in} In agreement with Taneja and Tuteja (1986), here we consider the weighted inaccuracy measure
\begin{eqnarray}\label{eq5}H_{X,Y}^w=-\int_0^\infty
x f(x)\ln g(x)dx,\end{eqnarray}
which is a quantitative-qualitative measure of inaccuracy associated with the statement of an experimenter. When $g(x) = f(x)$, then (\ref{eq5}) becomes (\ref{eq4}), the weighted entropy. For more properties of quantitative-qualitative measure of inaccuracy one may refer to Prakash and Taneja (1986) and Bhatia and Taneja (1991), among others. The following example illustrates the importance of qualitative characteristic of information as reflected in the definition of weighted inaccuracy measure.
\begin{ex}\label{ex1.1} Let $X_1$ and $Y_1$ denote random lifetimes of two
components with probability density functions
$f_1(x)=x/2,~x\in(0,2)$ and $g_1(x)=(2-x)/2,~x\in(0,2)$
respectively. By simple calculations, we have
$H_{X_1,Y_1}=H_{Y_1,X_1}=3/2$. But,
$$H_{X_1,Y_1}^w=\frac{22}{9}~{\rm and}~H_{Y_1,X_1}^w=\frac{5}{9}.$$
Therefore, the inaccuracy measure of the observer for the
observations $X_1$ (resp. $Y_1$) taking $Y_1$ (resp. $X_1$) as
corresponding assigned outcomes by the experimenter are identical.
Instead, $H_{X_1,Y_1}^w>H_{Y_1,X_1}^w$, i.e., weighted inaccuracy of
the observer for ($X_1,~Y_1$) is higher than that for
($Y_1,~X_1$). As a matter of fact, the inaccuracies measured from a quantitative point of view, neglecting the qualitative side, are identical. To distinguish them, we must take into account the qualitative characteristic as given in (\ref{eq5}). $\hfill\square$
\end{ex}
\hspace*{.2in} Analogous to weighted residual and past entropies Kumar et al. (2010) and Kumar and Taneja
(2012) introduced the notion of weighted residual inaccuracy
measure given by
\begin{eqnarray}\label{eq6}
H_{X,Y}^w(t)=-\int_t^\infty x\frac{f(x)}{\overline
F(t)}\ln\left(\frac{g(x)}{\overline G(t)}\right)dx
\end{eqnarray}
and weighted past inaccuracy measure given by
\begin{eqnarray}\label{eq7}
\overline
H_{X,Y}^w(t)=-\int_0^tx\frac{f(x)}{F(t)}\ln\left(\frac{g(x)}{G(t)}\right)dx,
\end{eqnarray}
and studied their properties in analogy with weighted residual
entropy and weighted past entropy, respectively. For $t=0$,
(\ref{eq6}) reduces to (\ref{eq5}) and for $t=\infty$, (\ref{eq7})
reduces to (\ref{eq5}). Various aspects of (\ref{eq6}) and (\ref{eq7})
have been discussed in Kundu (2014).\\

\hspace*{.2in} The rest of the paper is arranged as follows. In
Section 2 we introduce the concept of weighted interval inaccuracy
measure for doubly truncated random variables. We obtain upper and
lower bounds for weighted interval inaccuracy measure. In Section
3 we provide characterizations of quite a few useful continuous
distributions based on this newly introduced measure including its
uniqueness property. The effect of monotone transformations on the
weighted interval inaccuracy measure has been discussed in
Section~4.
\section{Weighted interval inaccuracy measure}
In the study of income distribution, the inequality is computed not
only for income greater/smaller than a fixed value but also for income
between two values. For example, in many practical situations, it
is of interest to study the inequality of a population eliminating high
(richest population) and low (poorest population) values, and therefore
doubly truncated populations are considered. In reliability theory and
survival analysis, often individuals
whose event time lies within a certain time interval are only
observed and one has information about the lifetime between two
time points. Thus, an individual whose event time is not in this
interval is not observed and therefore information on the subjects
outside this interval is not available to the investigator.
Accordingly, Kotlarski (1972) studied the conditional expectation
for the doubly truncated random variables. Later, Navarro and Ruiz
(1996) generalized the failure rate and the conditional
expectation to the doubly truncated random variables. For various
related results one may refer to Ruiz and Navarro (1996), Betensky
and Martin (2003), Sankaran and Sunoj (2004) among others.
Recently, Sunoj et al. (2009) and Misagh and Yari (2010, 2012)
studied the measure of uncertainty and conditional measure for
doubly truncated random variables and obtained some
characterization results. Furthermore, Misagh and Yari (2011)
explored the use of weighted information measures for doubly
truncated random variables. Motivated by this, we introduce the
notion of weighted
interval inaccuracy measure for doubly truncated random variables.\\
\hspace*{.2in} Let us consider two nonnegative absolutely
continuous doubly truncated random variables $(X|t_1\leqslant
X\leqslant t_2)$ and $(Y|t_1\leqslant Y\leqslant t_2)$ where
$(t_1,t_2)\in D=\{(u,v)\in\Re_+^2 :F(u)<F(v) ~{\rm and}~
G(u)<G(v)\}$. Then the interval inaccuracy measure of $X$ and $Y$
at interval $(t_1,t_2)$ is given by
\begin{eqnarray}\label{eq2.1}H_{X,Y}(t_1,t_2)&=&-\int_{t_1}^{t_2}
\frac{f(x)}{F(t_2)-F(t_1)}\ln\frac{g(x)}{G(t_2)-G(t_1)}dx.
\end{eqnarray}
When $g(x)=f(x)$, we obtain measure of uncertainty for doubly
truncated random variable as given in (2.6) and (2.7) of Sunoj et
al. (2009). Various aspects of interval inaccuracy measure have been discussed in
Kundu and Nanda (2014). To construct a shift-dependent dynamic
measure of inaccuracy, we use (\ref{eq2.1}) and define weighted
interval inaccuracy measure for two-sided truncated random
variables.
\begin{d1} The weighted interval inaccuracy measure of $X$ and $Y$ at interval $(t_1,t_2)$
is given by
\begin{eqnarray}\label{eq2.2}H_{X,Y}^w(t_1,t_2)=-\int_{t_1}^{t_2}
x\frac{f(x)}{F(t_2)-F(t_1)}\ln\frac{g(x)}{G(t_2)-G(t_1)}dx.
\end{eqnarray}
\end{d1}
\begin{r1}Clearly, $H_{X,Y}^w(0,t)=\overline H_{X,Y}^w(t)$,
$H_{X,Y}^w(t,\infty)=H_{X,Y}^w(t)$ and
$H_{X,Y}^w(0,\infty)=H_{X,Y}^w$ as given in (\ref{eq7}),
(\ref{eq6}) and (\ref{eq5}) respectively.$\hfill\square$
\end{r1}
\hspace*{.2in} The following example clarifies the effectiveness of the weighted interval inaccuracy measure.
\begin{ex} Let $X_1$, $Y_1$ be the random lifetimes as given in Example \ref{ex1.1}. Also let $X_2,~Y_2$ denote random lifetimes of two
components with probability density functions
$f_2(x)=2x,~x\in(0,1)$ and $g_2(x)=2(1-x),~x\in(0,1)$
respectively. Since $X_1$, $Y_1$ and $X_2,~Y_2$ belong to different domains, the use of weighted inaccuracy measure (\ref{eq5}) to compare them informatively is not interpretable. The weighted interval inaccuracy measure in the interval (0.2,0.8) are $H_{X_1,Y_1}^w(0.2,0.8)=-0.1143$ and $H_{X_2,Y_2}^w(0.2,0.8)=-0.2416$. Hence, the weighted interval inaccuracy measure between $X_1,~Y_1$ is greater than of it between $X_2,~Y_2$ in the interval (0.2,0.8).~$\hfill\square$
\end{ex}
An alternative way of writing (\ref{eq2.2}) is as follows:
\begin{eqnarray*}H_{X,Y}^w(t_1,t_2)=-\frac{1}{F(t_2)-F(t_1)}\int_{t_1}^{t_2}xf(x)\ln
g(x)dx+\frac{\ln\{G(t_2)-G(t_1)\}}{F(t_2)-F(t_1)}\int_{t_1}^{t_2}xf(x)dx,
\end{eqnarray*}
where the second integral on the right hand side is equal to
$$t_2F(t_2)-t_1F(t_1)-\int_{t_1}^{t_2}F(x)dx,~~{\rm or}~~ t_1\overline F(t_1)-t_2\overline F(t_2)+\int_{t_1}^{t_2}\overline F(x)dx.$$
The weighted interval inaccuracy measure can also be written as
\begin{eqnarray}\label{eq2.3}H_{X,Y}^w(t_1,t_2)&=&-\int_{t_1}^{t_2}\int_0^x
\frac{f(x)}{F(t_2)-F(t_1)}\ln\frac{g(x)}{G(t_2)-G(t_1)}dydx\nonumber\\
&=&t_1H_{X,Y}(t_1,t_2)+\int_{t_1}^{t_2}H_{X,Y}(x,t_2)dx.
\end{eqnarray}
Furthermore,
\begin{eqnarray}\label{eq2.4}H_{X,Y}^w(t_1,t_2)=t_2H_{X,Y}(t_1,t_2)-\int_{t_1}^{t_2}H_{X,Y}(t_1,y)dy,
\end{eqnarray}
where $H_{X,Y}(t_1,t_2)$ is the interval inaccuracy measure given
in (\ref{eq2.1}). Differentiating (\ref{eq2.3}) and (\ref{eq2.4})
with respect to $t_1$ and $t_2$, respectively, we obtain
$$\frac{\partial}{\partial t_1}H_{X,Y}^w(t_1,t_2)=t_1\frac{\partial}{\partial
t_1}H_{X,Y}(t_1,t_2)~~{\rm and}~~\frac{\partial}{\partial
t_2}H_{X,Y}^w(t_1,t_2)=t_2\frac{\partial}{\partial
t_2}H_{X,Y}(t_1,t_2).$$
\begin{r1} Weighted interval inaccuracy measure is increasing
(decreasing) in $t_1$ if and only if the interval inaccuracy
measure is increasing (decreasing) in $t_1$. The result also holds
for $t_2$.$\hfill\square$
\end{r1}
\hspace*{.2in} We decompose the weighted Kerridge inaccuracy
measure in terms of weighted residual, past and interval
inaccuracy measures on using the similar approach to that of
Misagh and Yari (2011).
\begin{r1}\label{r2.3} Let $X$ and $Y$ be two absolutely continuous nonnegative
random variables with $E(X)<\infty$. Then, for all
$0<t_1<t_2<\infty$, the weighted Kerridge inaccuracy measure can
be decomposed as \begin{eqnarray*}&&H_{X,Y}^w=F(t_1)\overline
H_{X,Y}^w(t_1)+\left[F(t_2)-F(t_1)\right]H_{X,Y}^w(t_1,t_2)+\overline F(t_2)H_{X,Y}^w(t_2)\\
&&-E(X)\left[F^*(t_1)\ln G(t_1)+\left\{F^*(t_2)-F^*(t_1)
\right\}\ln\left\{G(t_2)-G(t_1)\right\}+\overline
F^*(t_2)\ln\overline G(t_2)\right],\end{eqnarray*} which can be
interpreted as follows. The weighted inaccuracy measure can be
decomposed into four parts: $(i)$ the weighted inaccuracy measure
for random variables truncated above $t_1$, $(ii)$ the weighted
inaccuracy measure in the interval $(t_1,t_2)$ given that the item
has failed after $t_1$ but before $t_2$, $(iii)$ the weighted
inaccuracy measure for random variables truncated below $t_2$ and
$(iv)$ the pseudo inaccuracy for trivalent random variables which
determines whether the item has failed
before $t_1$ or in between $t_1$ and $t_2$ or after $t_2$.\\
When $t_1=t_2=t,$ then the above can be written as
$$H_{X,Y}^w=F(t)\overline H_{X,Y}^w(t)+\overline F(t)H_{X,Y}^w(t)
-E(X)\left[F^*(t)\ln G(t)+\overline F^*(t)\ln\overline
G(t)\right],$$ a result obtained by Kumar and Taneja
(2012).$\hfill\square$
\end{r1}
\hspace*{.2in} In virtue of Remark 2.2, below we obtain the bounds
for the interval inaccuracy measure based on the monotonic
behavior of the weighted interval inaccuracy measure. We first
give definitions of general failure rate (GFR), general
conditional mean (GCM) and geometric vitality function of a random
variable $X$ truncated at two points $t_1$ and $t_2$ where
$(t_1,t_2)\in D$. For details one may refer to Navarro and Ruiz
(1996), Nair and Rajesh (2000) and Sunoj et al. (2009).
\begin{d1}
The GFR functions of a doubly truncated random variable
$(X|t_1<X<t_2)$ are given by
$h^X_1(t_1,t_2)=\frac{f(t_1)}{F(t_2)-F(t_1)}$ and
$h^X_2(t_1,t_2)=\frac{f(t_2)}{F(t_2)-F(t_1)}$. Similarly
$h^Y_1(t_1,t_2)$ and $h^Y_2(t_1,t_2)$ are defined for the random
variable $(Y|t_1<Y<t_2)$.$\hfill\square$
\end{d1}
\begin{d1}The GCM of a doubly truncated random variable
$(X|t_1<X<t_2)$ is defined by
$$m_X(t_1,t_2)=E(X|t_1<X<t_2)=\frac{1}{F(t_2)-F(t_1)}\int_{t_1}^{t_2}xf(x)dx.$$
\end{d1}
\begin{d1}The geometric vitality function for doubly truncated random
variable $(X|t_1<X<t_2)$ is given by
$${\mathcal G}_X(t_1,t_2)=E\left(\ln X|t_1<X<t_2\right),$$
which gives the geometric mean life of $X$ truncated at two points
$t_1$ and $t_2$, provided $E(\ln X)$ is finite. The corresponding
weighted version of it is given by ${\mathcal
G}^w_X(t_1,t_2)=E\left(X\ln X|t_1<X<t_2\right).$$\hfill\square$
\end{d1}
When $H_{X,Y}^w(t_1,t_2)$ is increasing in each of the arguments keeping the other fixed,
then on differentiating (\ref{eq2.2}) with respect to $t_1$ and
$t_2$, we get
$$\frac{h_1^Y(t_1,t_2)}{h_1^X(t_1,t_2)}-\ln
h_1^Y(t_1,t_2)\leqslant H_{X,Y}(t_1,t_2)\leqslant
\frac{h_2^Y(t_1,t_2)}{h_2^X(t_1,t_2)}-\ln h_2^Y(t_1,t_2).$$ The
following proposition gives bounds for the weighted interval
inaccuracy measure. The proof follows from (\ref{eq2.2}) and hence
omitted.
\begin{p1}If $g(x)$ is decreasing in $x>0$, then $$-m_X(t_1,t_2)\ln
h_1^Y(t_1,t_2)\leqslant H_{X,Y}^w(t_1,t_2)\leqslant
-m_X(t_1,t_2)\ln h_2^Y(t_1,t_2).$$ For increasing $g(x)$ the above
inequalities are reversed.$\hfill\square$
\end{p1}

\hspace*{.2in} In the following two theorems we provide upper and
lower bounds for the weighted interval inaccuracy measure based on
monotonic behavior of the GFR functions of $Y$.
\begin{t1}\label{th2.1} For fixed $t_2$,\\
$(i)$ if $h_1^Y(t_1,t_2)$ is decreasing in $t_1$ then
$H_{X,Y}^w(t_1,t_2)\geqslant-m_X(t_1,t_2)\ln h_1^Y(t_1,t_2)$,\\
and $(ii)$ increasing $h_1^Y(t_1,t_2)$ in $t_1$ implies\\
$H_{X,Y}^w(t_1,t_2)\leqslant-m_X(t_1,t_2)\ln
h_1^Y(t_1,t_2)-\int_{t_1}^{t_2}\frac{xf(x)}{F(t_2)-F(t_1)}\ln\frac{G(t_2)-G(x)}{G(t_2)-G(t_1)}dx.$
\end{t1}
Proof: Note that (\ref{eq2.2}) can be written as
\begin{eqnarray}\label{eq2.5}H_{X,Y}^w(t_1,t_2)=-\int_{t_1}^{t_2}\frac{xf(x)\ln
h_1^Y(x,t_2)}{F(t_2)-F(t_1)}dx-\int_{t_1}^{t_2}\frac{xf(x)}{F(t_2)-F(t_1)}\ln\frac{G(t_2)-G(x)}{G(t_2)-G(t_1)}dx.\end{eqnarray}
$(i)$ For $t_1<x$,
$\ln\frac{G(t_2)-G(x)}{G(t_2)-G(t_1)}\leqslant0$ and $\ln
h_1^Y(x,t_2)\leqslant\ln h_1^Y(t_1,t_2)$ if $h_1^Y(t_1,t_2)$ is
decreasing in $t_1$. Then, from (\ref{eq2.5}), we obtain
\begin{eqnarray*}H_{X,Y}^w(t_1,t_2)&\geqslant&-\int_{t_1}^{t_2}\frac{xf(x)}{F(t_2)-F(t_1)}\ln
h_1^Y(x,t_2)dx\\&\geqslant&-m_X(t_1,t_2)\ln h_1^Y(t_1,t_2).
\end{eqnarray*}
$(ii)$ The second part follows easily from (\ref{eq2.5}) on using
the fact that $\ln h_1^Y(x,t_2)\geqslant\ln h_1^Y(t_1,t_2)$ for
$t_1<x$.$\hfill\square$
\begin{r1} In the above theorem if we take
$t_2=\infty$, then we get the lower (resp. upper) bound for the weighted
residual inaccuracy measure as obtained by Kumar et al. (2010) (resp. Kundu, 2014).$\hfill\square$
\end{r1}
\begin{ex}\label{ex2.2}Let $X$ be a nonnegative random variable with
probability density function
\begin{equation}\label{eq2.7}f(x)=\left\{\begin{array}{ll}2x,&
0<x<1,\\
 0,& otherwise\end{array}\right.
 \end{equation} and $Y$ is uniformly distributed over $(0,a)$. Then $m_X(t_1,t_2)=\frac{2(t_1^2+t_1t_2+t_2^2)}{3(t_1+t_2)}$,
$h_1^Y(t_1,t_2)=\frac{1}{(t_2-t_1)}$ and
$H_{X,Y}^w(t_1,t_2)=\frac{2(t_1^2+t_1t_2+t_2^2)\ln(t_2-t_1)}{3(t_1+t_2)}$.
Note that right hand side of part (ii) is $\geqslant
\frac{2(t_1^2+t_1t_2+t_2^2)\ln(t_2-t_1)}{3(t_1+t_2)}.$ It is
easily seen that part (ii) of the above theorem is fulfilled. For
part (i), let $X$ be uniformly distributed over $[\alpha,\beta]$ and let $Y$
follow Pareto-I distribution given by
\begin{eqnarray}\label{eq2.8}G(t)=1-\frac{\alpha}{t},~t>\alpha(>0).
\end{eqnarray} Then $m_X(t_1,t_2)=\frac{(t_1+t_2)}{2},~\alpha<t_1<t_2<\beta$ and $h_1^Y(t_1,t_2)=\frac{t_2}{t_1(t_2-t_1)},$ which is
decreasing in $t_1$, for fixed $t_2>2t_1$. Now
\begin{eqnarray*}H_{X,Y}^w(t_1,t_2)+m_X(t_1,t_2)\ln
h_1^Y(t_1,t_2)&=&2\left[\frac{1}{t_2-t_1}\int_{t_1}^{t_2}x\ln
xdx-\frac{(t_1+t_2)}{2}\ln t_1\right]\\&\geqslant&0,
\end{eqnarray*} and equality holds for $t_1\rightarrow t_2.$ Hence
part (i) is also fulfilled.$\hfill\square$
\end{ex}
\hspace*{.2in} The proof of the following theorem is analogous to Theorem \ref{th2.1} but for completeness we give a brief outline of
the proof.
\begin{t1} For fixed $t_1$, if $h_2^Y(t_1,t_2)$ is decreasing in $t_2$ then\\
$H_{X,Y}^w(t_1,t_2)\leqslant-m_X(t_1,t_2)\ln
h_2^Y(t_1,t_2)-\int_{t_1}^{t_2}\frac{xf(x)}{F(t_2)-F(t_1)}\ln\frac{G(x)-G(t_1)}{G(t_2)-G(t_1)}dx.$
\end{t1}
Proof: We write (\ref{eq2.2}) as
\begin{eqnarray}\label{eq2.6}H_{X,Y}^w(t_1,t_2)=-\int_{t_1}^{t_2}\frac{xf(x)\ln
h_2^Y(t_1,x)}{F(t_2)-F(t_1)}dx-\int_{t_1}^{t_2}\frac{xf(x)}{F(t_2)-F(t_1)}\ln\frac{G(x)-G(t_1)}{G(t_2)-G(t_1)}dx.\end{eqnarray}
Hence the result follows from (\ref{eq2.6}) on using the fact
that, for $x<t_2$, $\ln
h^Y_2(t_1,x)\geqslant\ln h^Y_2(t_1,t_2)$ when
$h_2^Y(t_1,t_2)$ is decreasing in $t_2$. $\hfill\square$
\begin{r1} If in the above theorem we take
$t_1=0$, then we get
$$\overline H_{X,Y}^w(t_2)\leqslant -\tau_F(t_2)\left[\ln\phi_G(t_2)+1\right]-\frac{G(t_2)}{F(t_2)}\int_0^{t_2}\frac{xf(x)}{G(x)}dx,$$
an upper bound to the weighted past inaccuracy measure as obtained
in Theorem 4.2 of Kumar and Taneja (2012).
\end{r1}
\begin{ex} Let $X$ be a nonnegative random variable with
probability density function as given in (\ref{eq2.7}) and let $Y$ be
uniformly distributed over $(0,a)$. Since
$h_1^Y(t_1,t_2)=h_2^Y(t_1,t_2)=\frac{1}{(t_2-t_1)}$, on using the
same argument as in Example \ref{ex2.2}, it can easily be shown that
the condition of the above theorem is fulfilled. $\hfill\square$
\end{ex}
\begin{r1} It is not difficult to see from (\ref{eq2.6}) that, for fixed $t_1$, if $h_2^Y(t_1,t_2)$ is increasing in $t_2$ then
$H_{X,Y}^w(t_1,t_2)\geqslant-m_X(t_1,t_2)\ln h_2^Y(t_1,t_2)$. But it also can be shown that for a random variable
with support $[0,\infty)$, $h_2^Y(t_1,t_2)$ may not be increasing in $t_2$. This condition can be achieved if either
the support of the random variable is $(-\infty,b]$ with $b>0$ or $[0,b]$ with $b<\infty$.
\end{r1}
\section{Characterizations based on weighted interval inaccuracy measure}
In the literature, the problem of characterizing probability
distributions has been investigated by many researchers. The
standard practice in modeling statistical data is either to derive
the appropriate model based on the physical properties of the
system or to choose a flexible family of distributions and then
find a member of the family that is appropriate to the data. In
both the situations it would be helpful if we find
characterization theorems that explain the distribution. In fact,
characterization approach is very appealing to both theoreticians
and applied workers. In this section we show that weighted
interval inaccuracy measure can uniquely determine the
distribution function. We also provide characterizations of quite
a few useful continuous distributions in terms of
weighted interval inaccuracy measure.\\
\hspace*{.2in} First we define the proportional hazard rate model
(PHRM) and proportional reversed hazard rate model (PRHRM). Let
$X$ and $Y$ be two random variables with hazard rate functions
$h_F(\cdot)$, $h_G(\cdot)$ and reversed hazard rate functions
$\phi_F(\cdot)$, $\phi_G(\cdot)$, respectively. Then $X$ and $Y$
are said to satisfy the PHRM (cf. Cox, 1959), if there exists
$\theta>0$ such that $h_G(t)=\theta h_F(t)$, or equivalently,
$\overline G(t)=\left[\overline F(t)\right]^\theta$, for some
$\theta$. This model has been widely used in analyzing survival
data; see, for instance, Cox (1972), Ebrahimi and Kirmani (1996),
Gupta and Han (2001) and Nair and Gupta (2007) among others.
Similarly, $X$ and $Y$ are said to satisfy PRHRM proposed by Gupta
et al. (1998) in contrast to the celebrated PHRM with
proportionality constant $\theta>0$, if
$\phi_G(t)=\theta\phi_F(t)$. Or, equivalently,
$G(t)=\left[F(t)\right]^\theta$, for some $\theta$. This model is
flexible enough to accommodate both monotonic as well as
non-monotonic failure rates even though the baseline failure rate
is monotonic. See Sengupta et al. (1999), Di Crescenzo (2000) or
Gupta and Gupta (2007) for some results on this model.\\
\hspace*{.2in} The general characterization problem is to obtain
when the weighted interval inaccuracy measure uniquely determines
the distribution function. We consider the following
characterization result. For characterization of a distribution by
using its GFR functions one may refer to Navarro and Ruiz (1996).
\begin{t1}For two absolutely continuous nonnegative random
variables $X$ and $Y$, when $H_{X,Y}^w(t_1,t_2)$ is increasing in
$t_1$ (for fixed $t_2$) and decreasing in $t_2$ (for fixed $t_1$)
and $h_i^Y(t_1,t_2)=\theta h_i^X(t_1,t_2),~\theta>0,~i=1,2$,
respectively, then $H_{X,Y}^w(t_1,t_2)$ uniquely determines
$F(x)$.
\end{t1}
Proof: Differentiating (\ref{eq2.2}) with respect to $t_i,~i=1,2$,
we have
$$\frac{\partial}{\partial
t_1}H_{X,Y}^w(t_1,t_2)=t_1h_1^X(t_1,t_2)\left[H_{X,Y}(t_1,t_2)+\ln\theta-\theta+\ln
h_1^X(t_1,t_2)\right]$$
$${\rm and,}~\frac{\partial}{\partial
t_2}H_{X,Y}^w(t_1,t_2)=-t_2h_2^X(t_1,t_2)\left[H_{X,Y}(t_1,t_2)+\ln\theta-\theta+\ln
h_2^X(t_1,t_2)\right].$$ Then for any fixed $t_1$ and arbitrary
$t_2$, $h_1^X(t_1,t_2)$ is a positive solution of the equation
$\eta(x_{t_2})=0$, where
$$\eta(x_{t_2})=t_1x_{t_2}\left[H_{X,Y}(t_1,t_2)+\ln\theta-\theta+\ln
x_{t_2}\right]-\frac{\partial}{\partial t_1}H_{X,Y}^w(t_1,t_2).$$
Similarly, for any fixed $t_2$ and arbitrary $t_1$,
$h_2^X(t_1,t_2)$ is a positive solution of the equation
$\zeta(y_{t_1})=0$, where
$$\zeta(y_{t_1})=t_2y_{t_1}\left[H_{X,Y}(t_1,t_2)+\ln\theta-\theta+\ln
y_{t_1}\right]+\frac{\partial}{\partial t_2}H_{X,Y}^w(t_1,t_2).$$
Differentiating $\eta(x_{t_2})$ and $\zeta(y_{t_1})$ with respect
to $x_{t_2}$ and $y_{t_1}$, respectively, we get
$\frac{\partial\eta(x_{t_2})}{\partial
x_{t_2}}=t_1\left[H_{X,Y}(t_1,t_2)+\ln\theta-\theta+1+\ln
x_{t_2}\right]$ and $\frac{\partial\zeta(y_{t_1})}{\partial
y_{t_1}}=t_2\left[H_{X,Y}(t_1,t_2)+\ln\theta-\theta+1+\ln
y_{t_1}\right].$ Furthermore, second order derivatives are
$\frac{\partial^2\eta(x_{t_2})}{\partial
x^2_{t_2}}=\frac{t_1}{x_{t_2}}>0$ and
$\frac{\partial^2\zeta(y_{t_1})}{\partial
y^2_{t_1}}=\frac{t_2}{y_{t_1}}>0.$ So, both the functions
$\eta(x_{t_2})$ and $\zeta(y_{t_1})$ are minimized at
$x_{t_2}=\exp\left[\theta-\ln\theta-1-H_{X,Y}(t_1,t_2)\right]=y_{t_1},$
respectively. Here $\eta(0)=-\frac{\partial}{\partial
t_1}H_{X,Y}^w(t_1,t_2)<0,$ since we assume that
$H_{X,Y}^w(t_1,t_2)$ is increasing in $t_1$, and also, when
$x_{t_2}\rightarrow\infty,~\eta(x_{t_2})\rightarrow\infty.$
Similarly $\zeta(0)=\frac{\partial}{\partial
t_2}H_{X,Y}^w(t_1,t_2)<0,$ and $\zeta(y_{t_1})\rightarrow\infty$
as $y_{t_1}\rightarrow\infty$. Therefore, both the equations
$\eta(x_{t_2})=0$ and $\zeta(x_{t_2})=0$ have unique positive
solutions $h_1^X(t_1,t_2)$ and $h_2^X(t_1,t_2)$, respectively.
Hence the proof is completed on using the fact that GFR functions
uniquely determine the distribution function (cf. Navarro and
Ruiz, 1996).$\hfill\square$\\

\hspace*{.2in} Now we provide characterization theorems for some
continuous distributions using GFR, GCM, geometric vitality
function and weighted interval inaccuracy measure under PHRM and
PRHRM. Below we characterize uniform distribution. Recall that
$\frac{\partial h_1^Y(t_1,t_2)}{\partial
t_2}=-h_1^Y(t_1,t_2)h_2^Y(t_1,t_2)$ and $\frac{\partial h_1^Y(t_1,t_2)}{\partial
t_1}=h_1^Y(t_1,t_2)\left(\frac{g'(t_1)}{g(t_1)}+h_1^Y(t_1,t_2)\right)$.
\begin{t1} Let $X$ and $Y$ be two absolutely continuous random variables
satisfying PRHRM with proportionality constant $\theta(>0)$. A
relationship of the form
\begin{eqnarray}\label{eq3.1}H_{X,Y}^w(t_1,t_2)+m_X(t_1,t_2)\ln
h_1^Y(t_1,t_2)=(1-\theta)\left[{\mathcal
G}^w_Z(t_1,t_2)-m_X(t_1,t_2)\ln(t_1-\alpha)\right],
\end{eqnarray} where ${\mathcal G}^w_Z(t_1,t_2)=E\left[X\ln(X-\alpha)|t_1<X<t_2\right]$
and $\alpha<t_1<t_2<\beta$, holds if and only if $X$ denotes the
random lifetime of a component with uniform distribution over
$(\alpha,\beta)$.
\end{t1}
Proof: The {\it if part} is obtained from (\ref{eq2.2}). To prove
the converse, let us assume that (\ref{eq3.1}) holds. Then from
definition we can write
\begin{eqnarray}\label{eq3.2}-\int_{t_1}^{t_2}xf(x)\ln\frac{g(x)}{G(t_2)-G(t_1)}dx
+\ln\frac{g(t_1)}{G(t_2)-G(t_1)}\int_{t_1}^{t_2}xf(x)dx\nonumber\\
=(1-\theta)\left[\int_{t_1}^{t_2}x\ln(x-\alpha)f(x)dx-\ln(t_1-\alpha)\int_{t_1}^{t_2}xf(x)dx\right].
\end{eqnarray}
Differentiating (\ref{eq3.2}) with respect to $t_i,~i=1,2$ we get,
after some algebraic calculations,
$$g(t_i)=k(t_i-\alpha)^{\theta-1},~i=1,2~{\rm and }~k>0~{\rm
(constant)},$$ or $g(t)=k(t-\alpha)^{\theta-1}$, which gives the
required result.$\hfill\square$
\begin{c1}Under PRHRM the relation $$H_{X,Y}^w(t_1,t_2)+m_X(t_1,t_2)\ln
h_2^Y(t_1,t_2)=(1-\theta)\left[{\mathcal
G}^w_Z(t_1,t_2)-m_X(t_1,t_2)\ln(t_2-\alpha)\right],$$ where
${\mathcal G}^w_Z(t_1,t_2)=E\left[X\ln(X-\alpha)|t_1<X<t_2\right]$
and $\alpha<t_1<t_2<\beta$ characterizes the uniform distribution
over $(\alpha,\beta)$. $\hfill\square$
\end{c1}
\hspace*{.2in} Next, we give a theorem which characterizes the
power distribution.
\begin{t1} For two absolutely continuous random variables $X$ and
$Y$ satisfying PRHRM with proportionality constant $\theta(>0)$, the relation
\begin{eqnarray}\label{eq3.3}H_{X,Y}^w(t_1,t_2)+m_X(t_1,t_2)\ln
h_1^Y(t_1,t_2)=(1-c\theta)\left[{\mathcal
G}^w_X(t_1,t_2)-m_X(t_1,t_2)\ln t_1\right],
\end{eqnarray}
for all $0<t_1<t_2<b$, characterizes the power distribution
\begin{equation}\label{eq3.4} F(t)=\left\{\begin{array}{ll}\left(\frac{t}{b}\right)^{c},&
0<t<b,\;b,c>0\\
 0,& otherwise.\end{array}\right.
 \end{equation}
\end{t1}
Proof: If $X$ follows power distribution as given in
(\ref{eq3.4}), then (\ref{eq3.3}) is obtained from (\ref{eq2.2}).
To prove the converse, let us assume that (\ref{eq3.3}) holds.
Then differentiating with respect to $t_i,~i=1,2$, we get, after
some algebraic calculations,
$$g(t_i)=kt_i^{c\theta-1},~i=1,2~{\rm and }~k>0~{\rm
(constant)},$$ or $g(t)=kt^{c\theta-1}$, which gives the required
result.$\hfill\square$
\begin{c1} The relationship $$H_{X,Y}^w(t_1,t_2)+m_X(t_1,t_2)\ln
h_2^Y(t_1,t_2)=(1-c\theta)\left[{\mathcal
G}^w_X(t_1,t_2)-m_X(t_1,t_2)\ln t_2\right]$$ characterizes the
power distribution as given in (\ref{eq3.4}) under PRHRM.
$\hfill\square$
\end{c1}
\hspace*{.2in} Below we characterize Weibull distribution under
PHRM.
\begin{t1}\label{th3.4}
Let $X$ and $Y$ be two absolutely continuous random variables
satisfying PHRM with proportionality constant $\theta(>0)$. A
relationship of the form
\begin{eqnarray}\label{eq3.5}H_{X,Y}^w(t_1,t_2)+m_X(t_1,t_2)\ln
h_1^Y(t_1,t_2)&=&(1-p)\left[{\mathcal
G}^w_X(t_1,t_2)-m_X(t_1,t_2)\ln
t_1\right]\nonumber\\&&+\lambda\theta\left[m_{X^{p+1}}(t_1,t_2)-t_1^pm_X(t_1,t_2)\right],
\end{eqnarray} where $m_{X^{p+1}}(t_1,t_2)=E\left(X^{p+1}|t_1<X<t_2\right)$, the conditional
expectation of $X^{p+1}$, holds for all $(t_1,t_2)\in D$ and $p>0$
if and only if $X$ follows Weibull distribution $$\overline
F(t)=e^{-\lambda t^p},~t>0,~p>0.$$
\end{t1}
Proof: The {\it if part} is straight forward. To prove the
converse, let us assume that (\ref{eq3.5}) holds. Then
differentiating with respect to $t_i,~i=1,2$, we get, after some
algebraic calculations,
$$g(t_i)=kt_i^{p-1}e^{-\lambda\theta t_i^p},~i=1,2~{\rm and }~k>0~{\rm
(constant)},$$ or $g(t)=kt^{p-1}e^{-\lambda\theta t^p}$, which
gives the required result.$\hfill\square$
\begin{c1} Under PHRM, the relation \begin{eqnarray*}H_{X,Y}^w(t_1,t_2)+m_X(t_1,t_2)\ln
h_2^Y(t_1,t_2)&=&(1-p)\left[{\mathcal
G}^w_X(t_1,t_2)-m_X(t_1,t_2)\ln
t_2\right]\\&&+\lambda\theta\left[m_{X^{p+1}}(t_1,t_2)-t_2^pm_X(t_1,t_2)\right],
\end{eqnarray*} where $m_{X^{p+1}}(t_1,t_2)=E\left(X^{p+1}|t_1<X<t_2\right)$, the conditional
expectation of $X^{p+1}$, characterizes the Weibull distribution
as given in the above theorem.$\hfill\square$
\end{c1}
\begin{r1} Taking $p=1$ in Theorem \ref{th3.4}, we obtain the
characterization theorem for exponential distribution with mean
$1/\lambda$. Similarly, $p=2$ characterizes the Rayleigh
distribution $\overline F(t)=e^{-\lambda t^2},
t>0$.$\hfill\square$
\end{r1}
\hspace*{.2in} Now we consider Pareto-type distributions which are
flexible parametric models and play important role in reliability,
actuarial science, economics, finance and telecommunications.
Arnold (1983) proposed a general version of this family of
distributions called Pareto-IV distribution having the cumulative
distribution function
\begin{eqnarray}\label{eq3.6}F(x)=1-\left[1+\left(\frac{x-\mu}{\beta}\right)^\frac{1}{\gamma}\right]^{-\alpha},~x>\mu,
\end{eqnarray}
where $-\infty<\mu<\infty,~\beta>0,~\gamma>0$ and $\alpha>0$.
This distribution is related to many other families of
distributions. For example, setting $\alpha=1,~\gamma=1$ and
$(\gamma=1,~\mu=\beta)$ in (\ref{eq3.6}), one at a time, we obtain
Pareto-III, Pareto-II and Preto-I distributions, respectively.
Also, taking $\mu=0$ and $\gamma\rightarrow\frac{1}{\gamma}$ in
(\ref{eq3.6}), we obtain Burr-XII distribution.\\
\hspace*{.2in} Now we consider Pareto-type distributions for
characterization under PHRM. Below we provide characterization of
 Pareto-I distribution.
 \begin{t1}\label{th3.5} Let $X$ and $Y$ be two absolutely continuous random variables
satisfying PHRM with proportionality constant $\theta(>0)$. Then the
relation
\begin{eqnarray}\label{eq3.7}H_{X,Y}^w(t_1,t_2)+m_X(t_1,t_2)\ln
h_1^Y(t_1,t_2)=(\alpha\theta+1)\left[{\mathcal
G}^w_X(t_1,t_2)-m_X(t_1,t_2)\ln t_1\right],
\end{eqnarray}
holds for all $\beta<t_1<t_2$ if and only if $X$ follows Pareto-I
distribution given by
\begin{eqnarray*}F(t)=1-\left(\frac{\beta}{t}\right)^\alpha,~t>\beta,~\alpha,\beta>0.
\end{eqnarray*}
\end{t1}
Proof: The {\it if part} is straightforward. To prove the
converse, let us assume that (\ref{eq3.7}) holds. Then
differentiating with respect to $t_i,~i=1,2$, we get, after some
algebraic calculations,
$$g(t_i)=kt_i^{-(\alpha\theta+1)},~i=1,2~{\rm and }~k>0~{\rm
(constant)},$$ or $g(t)=kt^{-(\alpha\theta+1)}$, which gives the
required result.$\hfill\square$
\begin{c1} The relation $$H_{X,Y}^w(t_1,t_2)+m_X(t_1,t_2)\ln
h_2^Y(t_1,t_2)=(\alpha\theta+1)\left[{\mathcal
G}^w_X(t_1,t_2)-m_X(t_1,t_2)\ln t_2\right]$$ characterizes the
same distribution under PHRM as mentioned in the above
theorem.$\hfill\square$
\end{c1}
\hspace*{.2in} We conclude this section by characterizing
Pareto-II distribution. The proof is similar to that of Theorem
\ref{th3.5} and hence omitted.
\begin{t1}
Let $X$ and $Y$ be two absolutely continuous random variables
satisfying PHRM with proportionality constant $\theta(>0)$. Then the
relation
\begin{eqnarray}H_{X,Y}^w(t_1,t_2)+m_X(t_1,t_2)\ln
h_1^Y(t_1,t_2)=(\alpha\theta+1)\left[{\mathcal
G}^w_Z(t_1,t_2)-m_X(t_1,t_2)\ln (t_1-\mu+\beta)\right],
\end{eqnarray}
where ${\mathcal
G}^w_Z(t_1,t_2)=E\left(X\ln(X-\mu+\beta)|t_1<X<t_2\right)$ holds
for all $\mu<t_1<t_2$ if and only if $X$ follows Pareto-II
distribution given by
\begin{eqnarray*}F(t)=1-\left[1+\left(\frac{t-\mu}{\beta}\right)\right]^{-\alpha},~t>\mu.
\end{eqnarray*}
\end{t1}
\begin{c1} Under PHRM the relation $$H_{X,Y}^w(t_1,t_2)+m_X(t_1,t_2)\ln
h_2^Y(t_1,t_2)=(\alpha\theta+1)\left[{\mathcal
G}^w_Z(t_1,t_2)-m_X(t_1,t_2)\ln (t_2-\mu+\beta)\right],$$ where
${\mathcal
G}^w_Z(t_1,t_2)=E\left[X\ln(X-\mu+\beta)|t_1<X<t_2\right]$ and
$\mu<t_1<t_2$ characterizes the same distribution as mentioned in
the above theorem.
\end{c1}
\section{Monotonic transformations}
In this section we study the weighted interval inaccuracy measure
under strict monotonic transformations. The following result is a
generalization of Theorem 4.1 of Di Cresenzo and Longobardi
(2006).
\begin{t1} Let $X$ and $Y$ be two absolutely continuous
nonnegative random variables. Suppose $\varphi(x)$ is strictly
monotonic, continuous and differentiable function with derivative
$\varphi'(x)$. Then, for all $0<t_1<t_2<\infty$,
\begin{equation*}H^w_{\varphi(X),\varphi(Y)}(t_1,t_2)=\left\{\begin{array}{ll}
H_{X,Y}^{w,\varphi}\left(\varphi^{-1}(t_1),\varphi^{-1}(t_2)\right)\\+E\left[\varphi(X)\ln\varphi'(X)|\varphi^{-1}(t_1)<X<\varphi^{-1}(t_2)\right],
~\varphi~strictly~increasing\\
H_{X,Y}^{w,\varphi}\left(\varphi^{-1}(t_2),\varphi^{-1}(t_1)\right)\\+E\left[\varphi(X)\ln\{-\varphi'(X)\}|\varphi^{-1}(t_2)<X<\varphi^{-1}(t_1)\right],
~\varphi~strictly~decreasing,\end{array}\right.
\end{equation*} where $$H_{X,Y}^{w,\varphi}(t_1,t_2)=-\int_{t_1}^{t_2}
\varphi(x)\frac{f(x)}{F(t_2)-F(t_1)}\ln\frac{g(x)}{G(t_2)-G(t_1)}dx.$$
\end{t1}
Proof: Let $\varphi(x)$ be strictly increasing. Then from
(\ref{eq5}), (\ref{eq6}) and (\ref{eq7}) we have
\begin{eqnarray}\label{eq4.1}H^w_{\varphi(X),\varphi(Y)}=H^{w,\varphi}_{X,Y}+E\left[\varphi(X)\ln\varphi'(X)\right],
\end{eqnarray}
where $H^{w,\varphi}_{X,Y}=-\int_0^\infty\varphi(x)f(x)\ln g(x)dx,$
\begin{eqnarray}\label{eq4.2}H^w_{\varphi(X),\varphi(Y)}(t)=H^{w,\varphi}_{X,Y}(\varphi^{-1}(t))+E\left[\varphi(X)\ln\varphi'(X)|X>\varphi^{-1}(t)\right],
\end{eqnarray}
where $H_{X,Y}^{w,\varphi}(t)=-\int_t^\infty
\varphi(x)\frac{f(x)}{\overline F(t)}\ln\left(\frac{g(x)}{\overline
G(t)}\right)dx$, and
\begin{eqnarray}\label{eq4.3}\overline H^w_{\varphi(X),\varphi(Y)}(t)=\overline H^{w,\varphi}_{X,Y}(\varphi^{-1}(t))+E\left[\varphi(X)\ln\varphi'(X)|X<\varphi^{-1}(t)\right],
\end{eqnarray}
where $\overline
H_{X,Y}^{w,\varphi}(t)=-\int_0^t\varphi(x)\frac{f(x)}{F(t)}\ln\left(\frac{g(x)}{G(t)}\right)dx$.
Now from Remark \ref{r2.3} we can write
\begin{eqnarray*}H^w_{\varphi(X),\varphi(Y)}&=&F(\varphi^{-1}(t_1))\overline
H^w_{\varphi(X),\varphi(Y)}(t_1)+\left[F(\varphi^{-1}(t_2))-F(\varphi^{-1}(t_1))\right]H^w_{\varphi(X),\varphi(Y)}(t_1,t_2)\\&&+\overline
F(\varphi^{-1}(t_2))H^w_{\varphi(X),\varphi(Y)}(t_2)-E(\varphi(X))\left[F^{w,\varphi}(\varphi^{-1}(t_1))\ln
G(\varphi^{-1}(t_1))\right.\\&&+\left.\{F^{w,\varphi}(\varphi^{-1}(t_2))-F^{w,\varphi}(\varphi^{-1}(t_1))\}
\ln\{G(\varphi^{-1}(t_2))-G(\varphi^{-1}(t_1))\}\right.\\&&+\left.\overline
F^{w,\varphi}(\varphi^{-1}(t_2))\ln\overline
G(\varphi^{-1}(t_2))\right],\end{eqnarray*} where
$F^{w,\varphi}(t)=\frac{1}{E[\varphi(X)]}\int_0^t\varphi(x)f(x)dx$. On using
(\ref{eq4.1}), (\ref{eq4.2}) and (\ref{eq4.3}) we obtain
\begin{eqnarray}\label{eq4.4}
&&H^{w,\varphi}_{X,Y}+E\left[\varphi(X)\ln\varphi'(X)\right]\nonumber=\left[F(\varphi^{-1}(t_2))-F(\varphi^{-1}(t_1))\right]H^w_{\varphi(X),\varphi(Y)}(t_1,t_2)\nonumber\\
&&+F(\varphi^{-1}(t_1))E\left[\varphi(X)\ln\varphi'(X)|X<\varphi^{-1}(t_1)\right]+\overline
F(\varphi^{-1}(t_2))E\left[\varphi(X)\ln\varphi'(X)|X>\varphi^{-1}(t_2)\right]\nonumber\\&&+F(\varphi^{-1}(t_1))\overline
H^{w,\varphi}_{X,Y}(\varphi^{-1}(t_1))+\overline
F(\varphi^{-1}(t_2))H^{w,\varphi}_{X,Y}(\varphi^{-1}(t_2))\nonumber\\
&&-E(\varphi(X))\left[F^{w,\varphi}(\varphi^{-1}(t_1))\ln
G(\varphi^{-1}(t_1))+\overline F^{w,\varphi}(\varphi^{-1}(t_2))\ln\overline
G(\varphi^{-1}(t_2))\right.\nonumber\\&&+\left.\{F^{w,\varphi}(\varphi^{-1}(t_2))-F^{w,\varphi}(\varphi^{-1}(t_1))\}
\ln\{G(\varphi^{-1}(t_2))-G(\varphi^{-1}(t_1))\}\right],
\end{eqnarray}
where the last three terms on the right hand side of (\ref{eq4.4})
are equal to
$$H^{w,\varphi}_{X,Y}-\left[F(\varphi^{-1}(t_2))-F(\varphi^{-1}(t_1))\right]H^{w,\varphi}_{X,Y}\left(\varphi^{-1}(t_1),\varphi^{-1}(t_2)\right),$$
giving the first part of the proof. If $\varphi(x)$ is strictly
decreasing we similarly obtain the second part of the
proof.$\hfill\square$
\begin{r1} Let $\varphi_1(x)=F(x)$ and $\varphi_2(x)=\overline F(x)$,
with $\varphi_1$ and $\varphi_2$ satisfying the assumptions of Theorem
4.1. Here $\varphi_1(X)$ and $\varphi_2(X)$ are uniformly distributed
over $(0,1)$. Then, for all $(t_1,t_2)\in D$, we have
$$H^w_{F(X),F(Y)}(t_1,t_2)=H_{X,Y}^{w,F}\left(F^{-1}(t_1),F^{-1}(t_2)\right)+E\left[F(X)\ln
f(X)|F^{-1}(t_1)<X<F^{-1}(t_2)\right]$$ and
$$H^w_{\overline F(X),\overline F(Y)}(t_1,t_2)=H_{X,Y}^{w,\overline F}\left(\overline F^{-1}(t_2),\overline F^{-1}(t_1)\right)
+E\left[\overline F(X)\ln f(X)|\overline F^{-1}(t_2)<X<\overline
F^{-1}(t_1)\right].$$
\end{r1}
\begin{r1} For two absolutely continuous nonnegative random variables $X$ and $Y$
$$H^w_{aX,aY}(t_1,t_2)=aH^w_{X,Y}\left(\frac{t_1}{a},\frac{t_2}{a}\right)+m_X\left(\frac{t_1}{a},\frac{t_2}{a}\right)a\ln
a$$ for all $a>0$ and $t_1>0$. Furthermore, for all $0<b<t_1$
$$H^w_{X+b,Y+b}(t_1,t_2)=H^w_{X,Y}(t_1-b,t_2-b)+bH_{X,Y}(t_1-b,t_2-b).$$
\end{r1}

\end{document}